\documentclass[a4paper,10pt]{scrartcl}

\usepackage[left=3.5cm, right=3.5cm, top=2.2cm, bottom=2.4cm]{geometry}
\usepackage[utf8]{inputenc}
\usepackage{babel}[german]
\usepackage[T1]{fontenc}
\usepackage{amsmath}
\usepackage{amssymb}
\usepackage{multicol}
\usepackage{graphicx}
\usepackage{xcolor}
\usepackage{lmodern}
\usepackage{enumerate}
\usepackage{caption}
\DeclareMathSizes{10}{9.5}{8}{7}

\title{Sufficient minimum degree conditions for the existence\\  of highly connected or edge-connected subgraphs}
\author{Maximilian Krone}
\date{May 2024}

\newcounter{ThNr}[section]
\DeclareRobustCommand{\newTh}[1]{\refstepcounter{ThNr}\arabic{ThNr} \label{#1}}

\begin{document}
\parskip2.5explus1exminus1ex
\parindent0mm

\begin{center}
\vspace*{6mm}
\textbf{\Large Sufficient minimum degree conditions for the existence\\  of highly connected or edge-connected subgraphs}

\vspace*{2ex}
{\large Maximilian Krone}

{\small Technische Universität Ilmenau\\ 
August 2026 %\today
} \vspace*{2ex}\end{center}

\begin{abstract}
Mader conjectured in 1979 that an average degree of at least $3k-1$ in a graph is sufficient for the existence of a $(k+1)$-connected subgraph. The following minimum degree analogue holds: Every graph with minimum degree at least $3k-1$ contains a $(k+1)$-connected subgraph on more than $2k$ vertices. Moreover, for triangle-free graphs, already an average degree of at least $2k$ is sufficient for a $(k+1)$-connected subgraph, which has at least $2(k+1)$ vertices. \vspace{1ex}

For edge-connectivity (in simple graphs), we prove the following: Every graph with average degree at least $2k$ contains a $(k+1)$-edge-connected subgraph on more than $2k$ vertices. Moreover, for every small $\alpha>0$ and for $k$ large enough in terms of $\alpha$, already a minimum degree of at least $k+k^{\frac{1}{2}+\alpha} = \big(1+o(1)\big)k$ is sufficient for a $(k+1)$-edge-connected subgraph.  \vspace{1ex}

It is shown that all of these results are sharp in some sense. The results are applied to decompose graphs into two highly connected or edge-connected parts.
\end{abstract} 

\hrulefill \vspace*{-1ex}

\subsection*{Introduction}

Let $G$ be a graph. We set $v(G) :=|V(G)|$ and $e(G):=|E(G)|$. We define its
\begin{itemize}
\item \textit{connectivity} \hfill $\kappa(G) \,:=\, \min\big\{ |S| : S\subseteq V(G),\, G-S \text{ is disconnected or } v(G-S)\leq 1\big\}$;
\item \textit{edge-connectivity} \hfill $\lambda(G) \,:=\, \min\big\{ |S| : S\subseteq E(G),\, G-S \text{ is disconnected or } v(G-S)\leq 1\big\}$;
\item \textit{minimum degree}  $\delta(G) := \min\big\{ d_G(v) : v\in V(G) \big\}$\, for $v(G)>0$, and $\delta(G):=0$ otherwise;
\item \textit{average degree} ~ $d(G) \,:=\, \frac{1}{v(G)} \sum_{v\in V(G)} d_G(v)$~ for $v(G)>0$, and $d(G):=0$ otherwise.
\end{itemize}
$G$ is $k$-\textit{connected} if $\kappa(G)\geq k$, and $k$-\textit{edge-connected} if $\lambda(G)\geq k$.
It is well-known and easy to verify that
$$\kappa(G) \,\leq\, \lambda(G) \,\leq\, \delta(G) \,\leq\, d(G)\,,$$
and the four graph metrics can attain values arbitrarily far apart from each other.

Mader first studied the following extremal question: which number of edges in a graph guarantees a $k$-connected \cite{Mader} or $k$-edge-connected \cite{Mader e} subgraph? Through $d(G)=\frac{2e(G)}{v(G)}$, every such condition on the edge number can be translated into some average degree condition: An average degree of at least $4k-2$ or $2k$ is sufficient for a $(k+1)$-connected or $(k+1)$-edge-connected subgraph, respectively.

Considering minimum degree conditions is a natural modification. Despite its title, the present paper also studies average degree conditions whenever even the respective minimum degree condition turns out to be best possible. Moreover, at some points, additional conditions on the size of the highly connected or edge-connected subgraph are investigated as they emerge naturally from the proofs.

Mader's result also implies that, if a graph contains a subgraph with high average degree, then it also contains a highly connected subgraph. The inverse statement is self-evident. Following this idea, both Mader's original results \cite{Mader, Mader e} and the main results from the present paper can be captured regarding the following four \textit{core metrics}, which were first investigated by Matula \cite{Matula}, seemingly independently of Mader's results:

We define the \textit{core connectivity} of a graph $G$ as $\kappa^*(G) \,:=\, \max\big\{ \kappa(H): H \text{ is a subgraph of G} \big\}$. Similarly, we define its \textit{core edge-connectivity} $\lambda^*(G)$, its \textit{core minimum degree} $\delta^*(G)$ (also known as \textit{degeneracy}), and its \textit{core average degree} $d^*(G)$ (also known as \textit{maximum average degree}).

\textbf{Theorem}\\
\textit{Let $G$ be a graph with $E(G)\neq\emptyset$.
\begin{enumerate}[(i)]
\item $1 \leq \kappa^*(G) \,\leq\, \lambda^*(G) \,\leq\, \delta^*(G) \,\leq\, d^*(G)$ \cite{Matula}.
%\item $1 \leq \hat{\kappa}(G) \,\leq\, \hat\lambda(G) \,\leq\, \hat\delta(G) \,\leq\, \hat d(G)$.
\item $d^*(G) \,<\, 4\kappa^*(G)-2$  ~\cite{Mader}.%[Korollar 1]
\item $\delta^*(G) \,<\, 3\kappa^*(G)-1$ ~(see Theorem \ref{Theorem}).
\item $d^*(G) \,<\, 2\lambda^*(G)$ \cite{Mader e, Matula}, ~and even the conclusion $d^*(G) \,<\, 2\delta^*(G)$ is sharp for every value of $\delta^*(G)$, even for bipartite graphs \cite{Mader e} (see also Theorem \ref{e Average}).
\item $\delta^*(G) \,=\, \big(1+o(1)\big)\lambda^*(G)$ ~for $\lambda^*(G)\to\infty$ ~ (see Theorem \ref{e Theorem}).
\item If $G$ is triangle-free, then $d^*(G) \,<\, 2\kappa^*(G)$ ~(see Theorem \ref{triangle-free}), and even the conclusion $\lambda^*(G) \,<\, 2\kappa^*(G)$ is sharp for every value of $\kappa^*(G)$, even for bipartite graphs.\\
\end{enumerate}}

Mader's conjecture in \cite{Mader Conj} is equivalent to $d^*(G) \,<\, 3\kappa^*(G)-1$, which would be sharp for every value of $\kappa^*(G)$ \cite{Mader, Matula}. Moreover, the author of the present paper conjectures that $\delta^*(G) \,<\, C\kappa^*(G)-1$ for some constant $C<3$. By Theorem \ref{Construction}(ii), a suitable constant must satisfy $C\geq 3-\frac{1}{6}$, which improves on the lower bound $\frac{7}{3}$ found by Matula \cite{Matula}. 

For each two of the three core metrics $\lambda^*(G), \delta^*(G), d^*(G)$, the asymptotically sharp bound on the factor between them is stated in (iv) and (v). Matula stated that the estimate $\delta^*(G) \,<\, 2\lambda^*(G)$, which is another consequence of (iv), is sharp for every value $\lambda^*(G)\leq 6$ \cite{Matula}. Considering this fact, statement (v) might be the most surprising of all of the above results, stating that $\delta^*(G)$ and $\lambda^*(G)$ align asymptotically.\\

In the last section, we apply all of the above results to a decomposition problem, which was introduced by Thomassen \cite{Thomassen}. The aim is to decompose a highly connected graph into two highly connected or edge-connected parts.

There exist some relations between the core metrics of a graph $G$ and other graph metrics, for example to its chromatic number $\chi(G)$:
It is well-known that a graph $G$ admits a greedy $(\delta^*(G)+1)$-coloring (of its vertex set), that is, $\chi(G) \leq \delta^*(G)+1$. It is even true that $\chi(G) \leq \lambda^*(G)+1$ \cite{Matula}.

\iffalse, which is a consequence of an exercise in Diestel \cite[Section 4]{Diestel}.

\textbf{Proposition\\}
\textit{Every graph without a $k$-coloring contains a $k$-edge-connected subgraph.}

\textbf{Proof.~} Let $G$ be a minimal counter-example. In particular, $G$ itself is not $k$-colorable and not $k$-edge-connected.
Hence, there is a partition of $V(G)$ into two sets $A$ and $B$ such that there are at most $k-1$ edges between $A$ and $B$. 

As both $G[A]$ and $G[B]$ are $k$-colorable and $G$ is not, there is some edge $ab$ with $a\in A$, $b\in B$. As $G$ is minimal, $G-ab$ has a $k$-coloring, that is, a decomposition of $V(G)$ into $k$ possibly empty anticliques $X_1,\dots,X_k$. We may assume $a,b\in X_1$. Since the combined number of edges from $A\cap X_1$ to $B\setminus X_1$ and from $B\cap X_1$ to $A\setminus X_1$ is at most $k-2$, there is some $i\in\{2,\dots,k\}$ such that both $X'_1 = (A\cap X_1) \cup (B\cap X_i)$ and $X'_i := (A\cap X_i) \cup (B\cap X_1)$ are anticliques. Hence, exchanging $X_1$ and $X_i$ with $X'_1$ and $X'_i$ yields a $k$-coloring of $G$, a contradiction. \hfill $\Box$\\ \fi

A consequence of Theorem \ref{Theorem} is $\chi(G)\leq \delta^*(G)+1 < 3\kappa^*(G)$, that is, every graph without a $(3k-1)$-coloring contains a $(k+1)$-connected subgraph. This result was recently shown by Nguyen in \cite{Nguyen}, as well as the addition that a slightly larger chromatic number is even sufficient for a $(k+1)$-connected subgraph that itself has a large chromatic number. 
Matula first gave a construction of graphs with $\chi(G) = 2\kappa^*(G)$ in \cite{Matula}, for each value of $\kappa^*(G)$, and asked whether the construction is extremal, that is, whether every graph without a $2k$-coloring contains a $(k+1)$-connected subgraph. This statement apparently holds for $k=1$ and was recently shown to be true for $k=2$ \cite{Bonnet}. 

While all questions on the extremal relation between the core edge-connectivity $\lambda^*(G)$ and the graph metrics $\delta^*(G), d^*(G), \chi(G)$ are answered now, the respective questions regarding core connectivity $\kappa^*(G)$ remain open, by now for almost half a century. \newpage

\subsection*{On highly connected subgraphs}

Mader conjectured in \cite{Mader Conj} that every graph with average degree at least $3k-1$ contains a $(k+1)$-connected subgraph. Proving this conjecture turned out to be quite difficult, as the best known sufficient linear factor was improved over many years, but never reached $3$. We consider a weaker version of the problem by requiring a minimum degree of at least $3k-1$.

In exchange for the stronger requirement, the $(k+1)$-connected subgraph is guaranteed to have a non-trivial size of more than $2k$ vertices, which emerges naturally from the proof. This bonus outcome cannot be expected under the weaker assumption of an average degree of $3k-1$, since by a construction by Carmesin \cite{Carmesin}, there are graphs with average degree $\big(3+\frac{1}{3}\big)k-\mathcal{O}(1)$ that do not contain a $(k+1)$-connected subgraph on more than $2k$ vertices. 

\textbf{Theorem \newTh{Theorem}}\\
\textit{Let $k\geq 1$ be fixed. Every graph with minimum degree at least $3k-1$ contains a $(k+1)$-connected subgraph on more than $2k$ vertices.}

\textbf{Proof.~} 
Let $G$ be a graph with $v(G)>2k$ without a $(k+1)$-connected subgraph on more than $2k$ vertices.

Let $v\in V(G)$ be a vertex of degree $d_G(v)$. We define a delimited penalty
$$m_G(v) ~:=~ \min\Big\{\max\big\{(3k-1)-d_G(v),\, 0\big\},\,k\Big\}~ \in \{0,\dots,k\}\,.$$
We prove by induction on $v(G)>2k$ that the total penalty satisfies $M(G) ~\geq~ 2k^2$, where
$$M(G) ~:=~ \sum_{v\in V(G)} m_G(v)\,.$$
Since $v(G)>2k$, $G$ itself is not $(k+1)$-connected. Hence it has a $k$-separation $(A,B)$, that is, two induced but not spanning subgraphs $A$ and $B$ with $V(A)\cup V(B) = V(G)$,\, $E(A)\cup E(B) = E(G)$ and $|V(A)\cap V(B)|=k$. By symmetry, we may assume that $v(A) \geq v(B)$.

Let $a\in V(A)\setminus V(B)$, $b \in V(B)\setminus V(A)$ and $s \in  V(A)\cap V(B)$. We have $d_G(a) = d_A(a)$ and hence $m_G(a)=m_A(a)$. Similarly, $m_G(b)=m_B(b)$.

\textbf{Case I:~} $v(A)\leq 2k$ and $v(B)\leq 2k$.

We have $d_G(s)\leq v(G)-1$, and hence $m_G(s) \geq \min\big\{(3k-1)-(v(G)-1),k\big\} = 3k-v(G)$.\linebreak It is $m_G(a) = m_G(b) = k$. Hence
$$M(G) ~\geq~ \sum_s \big(3k-v(G)\big) + \sum_a k  + \sum_b k ~=~ k(3k-v(G))+(v(G)-k)k ~=~ 2k^2\,.$$

\textbf{Case II:~} $v(A)>2k$ and $v(B)\leq 2k$.

We use the induction hypothesis for $A$. Since the function $x\mapsto \min\big\{\max\{3k-1\,-x,\, 0\},\,k\big\}$ is non-expanding, we have $m_A(s) - m_G(s) \leq \big|d_A(s)-d_G(s)\big| \leq |V(G)\setminus V(A)| = v(B)-k$. Again, $m_G(b)=k$. Hence
$$\begin{aligned}M(G) ~&\geq~ \sum_a m_A(a) + \sum_s \big(m_A(a)-(v(B)-k)\big) + \sum_b k \\
~&=~ M(A)-k(v(B)-k)+(v(B)-k)k ~=~ M(A) ~\geq~ 2k^2\,.\end{aligned}$$

\textbf{Case III:~} $v(A)>2k$ and $v(B)>2k$.

We use the induction hypothesis for $A$ and $B$. We have $m_G(s) \geq 0 \geq m_A(s)+ m_B(s)-2k$. Hence
$$\begin{aligned}M(G) ~&\geq~ \sum_a m_A(a) + \sum_s \big(m_A(s)+ m_B(s)-2k\big) + \sum_b m_B(b) \\
~&=~ M(A)+M(B)-k(2k) ~\geq~ 2k^2+2k^2-2k^2 ~=~ 2k^2\,.\end{aligned}$$

This finishes the proof: If $G$ has minimum degree at least $3k-1$, then $v(G)>2k$ and $M(G)=0$, a contradiction. \hfill $\Box$\\

Can we expect a lower sufficient minimum degree if we drop the condition on the size of the $(k+1)$-connected subgraph? Using ideas from the following proof one can show that the minimum degree condition from Theorem \ref{Theorem} is best possible for $k\leq 6$, that is, there are graphs with minimum degree $3k-2$ without a $(k+1)$-connected subgraph. This is not obvious but shall not be featured here. Instead we now consider general $k$. We prove that it does not matter (essentially) whether we require size $2k$ or only size $\frac{4}{3}k$ for the $(k+1)$-connected subgraph:

\textbf{Theorem \newTh{Construction}} \vspace*{-2ex}
\textit{\begin{enumerate}
%$\mathrm{(i)}$
\item[(i)] There exists a graph with minimum degree $3k-3$ without a $(k+1)$-connected subgraph on more than $\big\lceil\frac{4}{3}k\big\rceil$ vertices. \vspace*{-2ex}
\item[(ii)] There exists a graph with minimum degree $\big(3-\frac{1}{6}\big)k-\mathcal{O}(1)$ without a $(k+1)$-connected subgraph.
\end{enumerate}}

\textbf{Proof.~} We present a construction for (i), which we will modify for (ii) afterwards. 

The construction works in two steps. First, we construct a graph $H$ without a $(k+1)$-connected subgraph on more than $\big\lceil\frac{4}{3}k\big\rceil$ vertices in which there is a set $X$ of $2k-1$ vertices of degree at least $2k-2$ and all vertices in $V(H)\setminus X$ have degree at least $3k-3$.

Let $a,b,c \in \mathbb{N}$ which differ by at most $1$ with $a+b+c=k-1$. Let $V_0$ be an anticlique on $k$ vertices and let $V_1,\dots,V_l$ be disjoint cliques on $k-1$ vertices that are also disjoint to $V_0$, where $l= 3+2a+b$. We partition each $V_j$ for $j\geq 1$ into three cliques $V_j = A_j\cup B_j\cup C_j$ with $|A_j|=a$, $|B_j|=b$ and $|C_j|=c$.

Let $w_4,w_5,\dots,w_l$ be a list of the $l-3=2a+b$ vertices in $W:=A_1\cup B_1 \cup A_2$. For $j\in \{1,\dots,l\}$, we add all possible edges between $V_j$ and some set $W_j \subseteq V_0\cup\dots \cup V_{j-1}$ of $k$ vertices. We set $W_j = V_0$ for $j\in\{1,2,3\}$ and $W_j = A_{j-1}\,\cup\, B_{j-2} \,\cup\, C_{j-3} \,\cup\, \{w_j\}$ otherwise. For $j\in\{1,2,3,l\}$, these edges are shown in Figure 1.

Note that there are no edges between each two of the sets $A_{j-1}$, $B_{j-2}$ and $C_{j-3}$.
Let $H$ be the graph obtained this way and let $X := C_{l-2}\cup B_{l-1}\cup C_{l-1} \cup V_l$. We have that 
$$|X| ~=~ c+b+c+(k-1) ~\leq~ 1+a+b+c+(k-1) ~=~ 2k-1\,.\vspace{1.5ex}$$ 

\begin{center}
\vspace*{-1ex}\includegraphics[width = 0.97\textwidth]{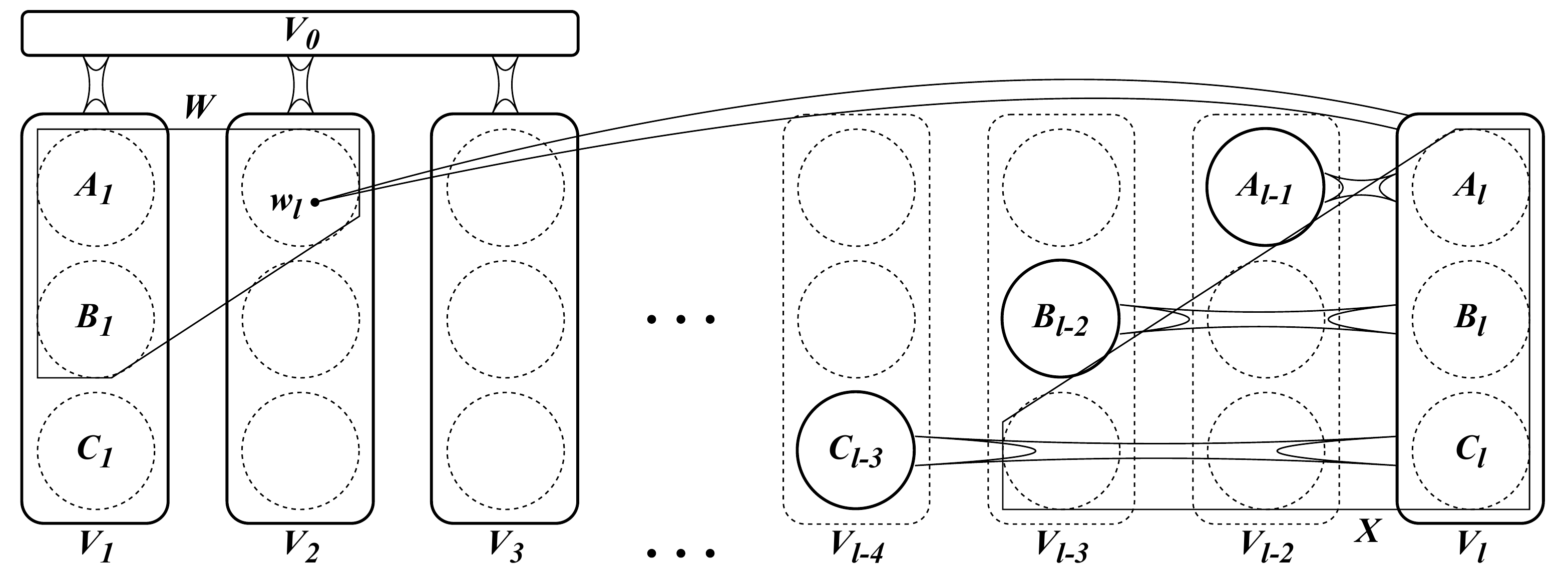}
\small Figure 1:~ The graph $H$. Only the edges between $V_j$ and $W_j$ for $j\in\{1,2,3,l\}$ are shown. \vspace{1.5ex}
\end{center}

We first check the claim on the degrees: All vertices in $V_0$ have the neighborhood $V_1\cup V_2 \cup V_3$ of size $3(k-1)$. The vertices in some $V_j$ for $j\geq 1$ have edges to the other $k-2$ vertices in $V_j$ and to the $k$ vertices in $W_j$. Hence, the vertices in $X$ have degree $2k-2$. Additionally, all vertices in $(V_1\cup\dots\cup V_l)\setminus X$ are contained in some set $W_j$, so they have $k-1$ additional neighbors in $V_j$, which adds up to a degree of $3k-3$.

We now prove inductively for $j\in \{3,4,\dots,l\}$ that there is no $(k+1)$-connected subgraph on more than $\big\lceil\frac{4}{3}k\big\rceil$ vertices in $G_j := H[V_0\cup\dots\cup V_j]$.

For $j=3$ this is true, since $G_3$ is separated by removing the $k$ vertices from $V_0$, so the vertex set of a $(k+1)$-connected subgraph must be completely contained in $V_0\cup V_1$, $V_0\cup V_2$ or $V_0\cup V_3$. But in each of these cases, the vertices in $V_0$ have degree at most $k$ and $V_1$, $V_2$ or $V_3$ by themselves are too small for containing a $(k+1)$-connected subgraph.

Now let $j\geq 4$ and suppose that there is a $(k+1)$-connected subgraph $G'$ on more than $\big\lceil\frac{4}{3}k\big\rceil$ vertices in $G_j$. Since the neighborhood $W_j$ of $V_j$ in $G_j$ consists of only $k$ vertices, $G'$ is either a subgraph of $G_j-V_j = G_{j-1}$ or of $G_j[V_j\cup W_j]$. The first case is prevented by the induction hypothesis, so suppose the second. Since there are no edges between each two of the sets $A_{j-1}$, $B_{j-2}$ and $C_{j-3}$, and those three are separated after deleting the set $V_j\cup \{w_j\}$ of $k$ vertices, $V(G')$ must be completely contained in one of the sets $V_j\cup \{w_j\} \cup A_{j-1}$,\, $V_j\cup \{w_j\} \cup B_{j-2}$\, or $V_j\cup \{w_j\} \cup C_{j-3}$. But all of them have size at most $(k-1)+1+\big\lceil\frac{k-1}{3}\big\rceil \leq \big\lceil\frac{4}{3}k\big\rceil$, a contradiction.

This proves the claimed properties of $H$. We now inductively construct a series $H_0,H_1,\dots$ of graphs without a $(k+1)$-connected subgraph on more than $\big\lceil\frac{4}{3}k\big\rceil$ vertices such that in $H_j$ there is a set $X_j$ of $\max\{2k-2^j, 0\}$ vertices of degree at least $2k-2$ and all vertices in $V(H_j)\setminus X_j$ have degree at least $3k-3$. Thus for $j\geq \log_2 2k$, $H_j$ has minimum degree at least $3k-3$.

Let $H_0 = H$ and $X_0=X$ (possibly replenished to size $2k-1$), which satisfy the conditions. Now suppose we have already constructed $H_j,X_j$ and we want to construct $H_{j+1},X_{j+1}$. Let $Y_j \subseteq X_j$ with $|Y_j| = \min\{k,|X_j|\}$. Let $H'_j$ be a copy of $H_j$. We define $H_{j+1}$ as the union of $H_j$ and $H'_j$ when we identify $Y_j$ and its copy in $H'_j$, that is, the subgraphs $H_j$ and $H'_j$ of $H_{j+1}$ intersect only in the vertex set $Y_j$. Since $|Y_i|\leq k$ and both $H_j$ and $H'_j$  contain no $(k+1)$-connected subgraph on more than $\big\lceil\frac{4}{3}k\big\rceil$ vertices, neither $H_{j+1}$ does. For every $y\in Y_j$, we have
$$d_{H_{j+1}}(y) ~=~ d_{H_j}(y) + d_{H'_j}(y) - d_{H_j[Y_j]}(y) ~\geq~ 2(k-1)+2(k-1)-(k-1) ~=~ 3(k-1)\,.$$
Also the vertices that are not contained in $X_j$ or its copy in $H_j'$ have degree at least $3k-3$. The remaining vertices from $X_j\setminus Y_j$ and its copy in $H_j'$ form the new set $X_{j+1}$ of size 
$$\begin{aligned} |X_{j+1}| ~&=~ 2|X_j\setminus Y_j| ~=~ 2\big(|X_j|-\min\{k, |X_j|\}\big) ~=~ 2\max \big\{|X_j|-k, 0\big\}\\
 ~&=~ 2\max \big\{\max \{2k-2^j, 0\} -k, 0\big\} ~=~ 2\max \big\{2k-2^j-k, 0\big\} ~=~ \max \{2k-2^{j+1}, 0\}\,. \end{aligned}$$
This finishes the proof of (i).

For (ii), we replace the complete graphs on each of the sets $A_j,B_j,C_j$ for $j\geq 1$ by the disjoint union of two complete graphs of size as equal as possible. By doing so, we lower the degree of every vertex in $H_0=H$ by at most $\frac{1}{6}k +\mathcal{O}(1)$. One observes inductively for $j=0,1,\dots$ that $d_{H_j[Y_j]}(v)\leq d_{H_j[X_j]}(v)\leq k-\frac{1}{6}k +\mathcal{O}(1)$ for every vertex $v\in Y_j$. Hence indeed  %For the construction of $H_1$, we may assume $c\leq a$, and we let $Y_0=V_l\subseteq X$. Hence, $|X_1|$ consists of $C_{l-2}\cup B_{l-1}\cup C_{l-1}$ and its copy, which are $2(2c+b)\leq 2(a+b+c) = 2k-2$ vertices. Doing so guarantees  for every $j\geq 0$, which is needed to show for $y\in Y_j$ that
$$\begin{aligned} d_{H_{j+1}}(v) ~&=~ d_{H_j}(v) + d_{H'_j}(v) - d_{H_j[Y_j]}(v)\\
~&\geq~ 2\big(2k-\tfrac{1}{6}k-\mathcal{O}(1)\big)-\big(k-\tfrac{1}{6}k +\mathcal{O}(1)\big) ~=~ \big(3-\tfrac{1}{6}\big)k-\mathcal{O}(1)\,.\end{aligned}$$
The proof of the non-existence of a $(k+1)$-connected subgraph first works as for (i): The vertex set of a $(k+1)$-connected subgraph $G'$ must be completely contained in one of the sets $V_j\cup \{w_j\} \cup A_{j-1}$,\, $V_j\cup \{w_j\} \cup B_{j-2}$\, or $V_j\cup \{w_j\} \cup C_{j-3}$, for some $j\geq 4$. Suppose $V(G') \subseteq V_j\cup \{w_j\} \cup A_{j-1}$. The two components of $H[A_{j-1}]$ can be separated by deleting the set $V_j\cup \{w_j\}$ of $k$ vertices. Hence $G'$ cannot contain vertices from both components. The same holds for $A_j \subseteq V_j$. Together, $V(G') \leq 2\big\lceil\frac{a}{2}\big\rceil+b+c+1 \leq 1+a+b+c+1 = k+1$, which is too low for being $(k+1)$-connected. The other two cases work similarly. \hfill$\Box$

One can check that the construction even produces a graph with edge-connectivity $3k-3$ or $\big(3-\frac{1}{6}\big)k-\mathcal{O}(1)$, respectively.

\textbf{Conjecture\\}
\textit{There is a $C<3$ such that, for all $k\in \mathbb{N}$, every graph with minimum degree at least $Ck-1$ contains a $(k+1)$-connected subgraph.}\vspace*{1ex}

By Theorem \ref{Construction}(ii), a suitable constant must satisfy $C\geq 3-\frac{1}{6}$. The conjecture is posed in a way that allows to assume that $k \geq k_0$ for some sufficiently large $k_0$: For $k < k_0$, $\big\lceil\big(3-\frac{1}{k_0}\big)k -1\big\rceil \,=\, \big\lceil 3k-\frac{k}{k_0} -1 \big\rceil \,=\, 3k-1$, so Theorem \ref{Theorem} suffices for $C=3-\frac{1}{k_0} < 3$.\\

As another modification of Mader's Conjecture, we now consider triangle-free graphs.

\textbf{Theorem \newTh{triangle-free}\\}
\textit{Let $k\geq 1$. Every triangle-free graph with average degree at least $2k$ contains a $(k+1)$-connected subgraph (on at least $2(k+1)$ vertices).}

\textbf{Proof.~} First, we prove the well-known statement that, if $G$ is a triangle-free graph on $2l$ vertices, then $e(G)\leq l^2$. This is true for $l=1$. If $l>1$ and $E(G)\neq \emptyset$, consider an edge $vw\in E(G)$. As $G$ is triangle-free, the neighborhoods of $v$ and $w$ are disjoint and hence $d_G(v)+d_G(w)\leq v(G)=2l$. This yields $e(G) = e(G-v-w)+d_G(v)+d_G(w)-1\linebreak \leq (l-1)^2+2l-1 = l^2$.

Now, let $G$ be a triangle-free graph with $v(G)\geq 2k$ without a $(k+1)$-connected subgraph. We prove by induction on $v(G)$ that $e(G) \leq k\big(v(G)-k\big)$. By the consideration above, the claim holds for $v(G)=2k$.

Now let $v(G)>2k$, for which $G$ itself is not $(k+1)$-connected. Hence it has a $k$-separation $(A,B)$ (see Theorem \ref{Theorem}). By symmetry, we may assume that $v(A) \geq v(B)$.

\textbf{Case I:~} $v(B) \leq 2k$.

Suppose that all vertices in $V(B)\setminus V(A)$ have degree greater than $k$. Let $b\in V(B)\setminus V(A)$. Since $|V(A)\cap V(B)|=k$, $b$ has a neighbor $c\in V(B)\setminus V(A)$. As $B$ is triangle-free, the neighborhoods of $b$ and $c$ are disjoint subsets of $V(B)$. By the assumption, both contain at least $k+1$ vertices, which contradicts $v(B) \leq 2k$.

Hence, there is a vertex $b\in V(B)\setminus V(A)$ of degree at most $k$. By using the induction hypothesis for $G-b$,
$$e(G) ~=~ e(G-b)+d_G(b) ~\leq~ e(G-b)+k ~\leq~ k\big(v(G)-1-k\big)+k ~=~ k\big(v(G)-k\big)\,.$$
\textbf{Case II:~} $v(B) > 2k$.

Thus also $v(A) > 2k$, so we may use the induction hypothesis for both $A$ and $B$. Hence,
$$e(G) ~\leq~ e(A)+e(B) ~\leq~ k\big(v(A)-k\big)+k\big(v(B)-k\big) ~=~ k\big(v(G)-k\big)\,.$$
This proves the claim. 

Hence, if $G$ is a triangle-free graph with average degree at least $2k$, then $v(G)\geq 2k$ and $e(G) \geq kv(G) > k\big(v(G)-k\big)$. Thus $G$ contains a $(k+1)$-connected subgraph $G'$.
As $G'$ is also triangle-free, for some edge $vw\in E(G')$, the neighborhoods of $v$ and $w$ are disjoint. This shows $v(G')\geq d_{G'}(v)+d_{G'}(w) \geq 2(k+1)$. \hfill $\Box$\\

Already in \cite{Mader}, Mader proved the above result for graphs with girth larger than $k+1$. He also presented a construction of bipartite graphs with arbitrarily large girth and average degree arbitrarily close to $2k$ without a $(k+1)$-connected subgraph. 

Theorem \ref{triangle-free} is still best possible when we replace the average degree by the minimum degree or even the edge-connectivity, that is, there exists a bipartite, $(2k-1)$-edge-connected graph without a $(k+1)$-connected subgraph. One can build such a graph by gluing together copies of $H_0 := K_{k,2k-1}$ similarly to the final step of the proof of Theorem \ref{Construction}. 

\subsection*{\boldmath On highly edge-connected subgraphs}

Mader proved in \cite{Mader e} that every graph $G$ with $e(G) > k\big(v(G)-\frac{k+1}{2}\big)$ and $v(G)\geq k$ contains a $(k+1)$-edge-connected subgraph. This condition is sharp for every $k$ and every $v(G)$: Consider the graph $G$ on a clique $A$ on $k$ vertices and an anticlique $V(G)\setminus A$ with all possible edges between $A$ and $V(G)\setminus A$. As all vertices in $V(G)\setminus A$ have degree $k$, $G$ does not even contain a subgraph with minimum degree at least $k+1$.
Note that, when $v(G)$ tends to $\infty$, the average degree of $G$ tends to $2k$. This also happens for the complete bipartite graph that is obtained by removing all edges from the clique $A$.

We first show that an average degree of at least $2k$ even implies the existence of a $(k+1)$-edge-connected subgraph on more than $2k$ vertices.

\textbf{Theorem \newTh{e Average}\\}
\textit{If $G$ is a loopless multigraph with $e(G)\,>\,k\,(v(G)-1)$, then $G$ contains a submultigraph $G'$ with \vspace*{-1ex}
\begin{enumerate}[(i)]
	\item $e(G')\,>\,k\,(v(G')-1)$. \vspace*{-1ex}
	\item $G'$ is $(k+1)$-edge-connected. \vspace*{-1ex}
	\item $G'$ contains $k$ pairwise edge-disjoint spanning trees. \vspace*{-1ex}
	\item If $G$ is a graph, then $v(G')> 2k$. \vspace*{-1ex}
	\item If $G$ is a graph, then $G'$ is $2$-connected. \vspace*{-1ex}
	\item If $G$ is a triangle-free graph, then $v(G')\geq 4k$.
\end{enumerate}}

\textbf{Proof.~} \begin{enumerate}[(i)]
	\item Let $G'$ be a smallest nonempty induced submultigraph of $G$ that satisfies the property $e(G')\,>\,k\,(v(G')-1)$.
	\item Let $\mathcal{P}$ be a partition of $V(G')$ in more than one class. Thus for all $A \in \mathcal{P}$, we have $e(G[A]) \leq k(|A|-1)$. Hence the number of edges between distinct classes of $\mathcal{P}$ satisfies
$$\begin{aligned} e(\mathcal{P}) ~&=~ e(G')\,-\sum_{A\in \mathcal{P}} e(G[A]) ~>~ k(v(G')-1)\,- \sum_{A\in \mathcal{P}} k(|A|-1)\\
~&=~ k\Big(v(G')-1+|\mathcal{P}|-\sum_{A\in \mathcal{P}} |A| \Big) ~=~ k(|\mathcal{P}|-1)\,. \end{aligned}$$
In particular, it follows from the case $|\mathcal{P}|=2$, that $G'$ is $(k+1)$-edge-connected.
	\item Moreover, by the Theorem by Tutte \cite{Tutte} and Nash-Williams \cite{N-W}, $G'$ contains $k$ pairwise edge-disjoint spanning trees.
	\item If $G$ and hence $G'$ is a (simple) graph, then $G'$ has more than $2k$ vertices, since otherwise $e(G') \,\leq\, \frac{1}{2}v(G')(v(G')-1) \,\leq\, k\,(v(G')-1)$, a contradiction. 
	\item Suppose that $G'$ is not $2$-connected. Since $v(G')\geq 3$, there is a vertex $s$ and a partition of $V(G')\setminus \{s\}$ into two sets $A$ and $B$ such that there is no edge between $A$ and $B$ in $G'$. Hence
	$$e(G') ~=~ e\big(G[A\cup\{s\}]\big) + e\big(G[B\cup\{s\}]\big) ~\leq~ k|A| + k|B| ~=~ k\,(v(G')-1)\,,$$
	a contradiction.
	\item If $G$ and hence $G'$ is additionally triangle-free, then it is well-known (see also the proof of Theorem \ref{triangle-free}) that
$$e(G') \leq \left\{ \begin{matrix}
\big(\frac{v(G')}{2}\big)^2 \leq \frac{v(G')+2}{2}\,\frac{v(G')-1}{2} & \text{ for even } v(G')\,,\vspace*{1ex}\\
\frac{v(G')+1}{2}\,\frac{v(G')-1}{2} & \text{ for odd } v(G')\,.
\end{matrix} \right.$$
This yields $v(G')\geq 4k$, since otherwise we could bound the first factor in both cases by $2k$ and obtain $e(G') \leq 2k \frac{v(G')-1}{2} = k(v(G')-1)$, a contradiction. \hfill $\Box$ \\ 
\end{enumerate}

In particular, Theorem \ref{e Average} applies to graphs with minimum degree at least $2k$. Even the minimum degree version is best possible in some senses: 
\begin{enumerate}
\item[(1)] The complete graph on $2k$ vertices has minimum degree $2k-1$ but no $(k+1)$-edge-connected subgraph on more than $2k$ vertices. Similarly, the complete bipartite graph on two classes of size $2k-1$ has minimum degree $2k-1$ but no $(k+1)$-edge-connected subgraph on at least $4k$ vertices.
\item[(2)] If a graph $G$ on more than one vertex contains $k$ pairwise edge-disjoint spanning trees, then $e(G)\geq k(v(G)-1)$. Similarly to the proof of (iv), we can conclude $v(G)\geq 2k$ with equality if and only if $G$ is the complete graph on $2k$ vertices. Now, if $v(G)>2k$, then the average degree of $G$ is greater than $2k-1$:
$$d(G)\,=\,\frac{2e(G)}{v(G)} \,\geq\, \frac{2k(v(G)-1)}{v(G)} \,>\, \frac{2k(2k-1)}{2k} \,=\, 2k-1\,.$$
Every $(2k-1)$-regular graph without a complete graph on $2k$ vertices as a component does not contain a subgraph with average degree more than $2k-1$, and thus no subgraph with $k$ pairwise edge-disjoint spanning trees.
\item[(3)] There exists a multigraph with minimum degree $2k-1$ without a $(k+1)$-edge-connected submultigraph: Consider the multigraph $G$ on the vertex set $\{1,\dots,2(k+1)\}$ with $k-1$ parallel edges between each two consecutive vertices $i$ and $i+1$, one additional edge between $1$ and each of the $k$ vertices $\{2,\dots,k+1\}$, and between $2(k+1)$ and each of the $k$ vertices $\{k+2,\dots,2k+1\}$. Clearly, all vertices have degree $2k-1$. Between the two vertex sets $\{1,\dots,k+1\}$ and $\{k+2,\dots,2k+2\}$ there are only $k-1$ edges. Hence, if $G'$ is a $(k+1)$-edge-connected submultigraph of $G$, then either $V(G')\subseteq \{1,\dots,k+1\}$ or $V(G')\subseteq \{k+2,\dots,2k+2\}$. In both cases, there is a vertex of degree at most $k$, $\max V(G')$ or $\min V(G')$, respectively. This is a contradiction.
\item[(4)] A similar construction can be used for (simple) graphs. Together with ideas from Theorem \ref{e Construction}, one can show that for $k\leq 6$, requiring a minimum degree of $2k$ is still best possible: One can build graphs with minimum degree $2k-1$ without a $(k+1)$-edge-connected subgraph. A suitable construction is also sketched in \cite{Matula}.
For $k\leq 4$, one can even find bipartite examples. As both constructions allow no generalization, this will not be featured here. Instead, we now consider large $k$.\\
\end{enumerate} 

As the construction from the beginning of this section indicates, the most determining condition for a $(k+1)$-connected subgraph is that it has no vertex of degree at most $k$. We now prove that, for $k\to \infty$, already a minimum degree $\big(1+o(1)\big)k$ in a graph is sufficient for a $(k+1)$-edge-connected subgraph.

Let $\alpha \in \big(0,\frac{1}{2}\big)$ be fixed. All of the following definitions will depend on $k$, which will not be reflected in the notations. Let $m = m(k) \in [0,k]$ not necessarily an integer, which we will choose later. 

Let $G$ be a graph and $v \in V(G)$ of degree $d_G(v)$. Similarly to the proof of Theorem \ref{Theorem}, we define a delimited penalty
$$m_G(v) ~:=~ \min\Big\{\max\big\{k+m-d_G(v),\, 0\big\},\,m\Big\}~ \in [0,m]\,.$$
We consider the following weighted penalty of the graph $G$:
$$M(G) ~:=~ \max\Big\{ \sum_{j=1}^{v(G)} \alpha j^{\alpha-1}m_G(v_j) ~:~ v_1,\dots,v_{v(G)} \text{ is a list of all vertices of } G \Big\}\,.$$
The sum of the first $l$ weights can be upper-bounded as follows:
$$\sum_{j=1}^l \alpha j^{\alpha-1} ~=~ \sum_{j=1}^l \int_{j-1}^j \alpha j^{\alpha-1} dx ~\leq~ \sum_{j=1}^l \int_{j-1}^j \alpha x^{\alpha-1} dx ~=~ \int_0^l \alpha x^{\alpha-1} dx ~=~ l^\alpha\,.$$
Similarly, we obtain as a lower bound:
$$\sum_{j=1}^l \alpha j^{\alpha-1} ~=~ \sum_{j=1}^l \int_{j}^{j+1} \alpha j^{\alpha-1} dx ~\geq~ \sum_{j=1}^l \int_{j}^{j+1} \alpha x^{\alpha-1} dx ~=~ \int_1^{l+1} \alpha x^{\alpha-1} dx ~\geq~ l^\alpha-1\,.$$

\textbf{Lemma \ref{e Theorem}a}\\
\textit{If $G$ is a graph that contains a vertex $w$ of degree at most $k$, then}
$$M(G) ~\geq~ M(G-w)-k^\alpha\,.$$

\textbf{Proof.~} Let $v_1,\dots,v_{v(G)-1}$ be a list of $V(G-w) = V(G) \setminus \{w\}$ that attains \vspace*{-1ex}
$$M(G-w) ~=~ \sum_{j=1}^{v(G)-1} \alpha j^{\alpha-1}m_{G-w}(v_j)\,.$$
For every $j\in \{1,\dots,v(G)-1\}$, we have $d_{G}(v_j)-d_{G-w}(v_j) \in \{0,1\}$ as $G$ is simple. Since the function $x\mapsto \min\big\{\max\{k+m-x,\, 0\},\,m\big\}$ is monotone decreasing and non-expanding, we have $0 \leq m_{G-w}(v_j) -m_{G} (v_j) \leq d_{G}(v_j)-d_{G-w}(v_j)$. Hence $\Delta_j := m_{G-w}(v_j) -m_{G}(v_j) \in [0,1]$ and
$$\sum_{j=1}^{v(G)-1} \Delta_j ~\leq~ \sum_{j=1}^{v(G)-1} d_{G}(v_j)-d_{G-w}(v_j) ~=~ d_G(w) ~\leq~ k\,.$$
Consider the list $v_1,\dots,v_{v(G)-1},w$ of the vertices of $G$. This yields
$$\begin{aligned} M(G) ~&\geq~ \sum_{j=1}^{v(G)-1} \alpha j^{\alpha-1}m_G(v_j) + \alpha v(G)^{\alpha-1}m_G(w)\\
~&>~ \sum_{j=1}^{v(G)-1} \alpha j^{\alpha-1}\big(m_{G-w}(v_j)-\Delta_j\big) ~=~ M(G-w)-\sum_{j=1}^{v(G)-1} \alpha j^{\alpha-1}\Delta_j\,.\end{aligned}$$
Since $\sum_{j=1}^{v(G)-1} \Delta_j \leq k$ and the weights $\alpha j^{\alpha-1}$ are decreasing, the weighted sum over the $\Delta_j \in [0,1]$ is maximized if $\Delta_j = 1$ for $j \leq k$ and $\Delta_j = 0$ otherwise. That is, \vspace*{-1ex}
$$M(G) ~\geq~ M(G-w)-\sum_{j=1}^{k} \alpha j^{\alpha-1} ~\geq~ M(G-w)-k^\alpha\,,$$
using the upper bound on the sum of the first $k$ weights. \hfill $\Box$ \\

\textbf{Lemma \ref{e Theorem}b}\\
\textit{If $G$ is a graph that contains a set $S$ of at most $k$ edges such that $G-S$ is disconnected, so it decomposes into non-empty subgraphs $A$ and $B$, then}
$$M(G) ~\geq~ 2^{\alpha-1}\big(M(A)+M(B)\big) - m\Big(\frac{4k}{m}\Big)^\alpha\,.$$

\textbf{Proof.~} By saying that $G-S$ decomposes into $A$ and $B$ we mean that $V(G)$ is partitioned into $V(A)\cup V(B)$ and $E(G)\setminus S$ is partitioned into $E(A)\cup E(B)$. We do not require $A$ and $B$ to be connected.

Let $a_1,\dots,a_{v(A)}$ and $b_1,\dots,b_{v(B)}$ be lists of $V(A)$ and $V(B)$ that attain \vspace*{-1ex}
$$M(A) ~=~ \sum_{j=1}^{v(A)} \alpha j^{\alpha-1}m_A(a_j) ~~\text{ and }~~ M(B) ~=~ \sum_{j=1}^{v(B)} \alpha j^{\alpha-1}m_B(b_j)\,.$$
Clearly $d_A(a_j) = d_{G-S}(a_j)$ and hence $m_A(a_j) = m_{G-S}(a_j)$. Similarly $m_B(b_j) = m_{G-S}(b_j)$. We may assume by symmetry that $v(A) \leq v(B)$. By considering the list
$$a_1,b_1,a_2,b_2,\dots,a_{v(A)},b_{v(A)},b_{v(A)+1},\dots,b_{v(B)}$$
of all vertices of $V(G-S)=V(G)$, we obtain \vspace*{-1ex}
$$\begin{aligned} M(G-S) ~&\geq~ \sum_{j=1}^{v(A)} \alpha (2j-1)^{\alpha-1}m_A(a_j) + \sum_{j=1}^{v(A)} \alpha (2j)^{\alpha-1}m_B(b_j) + \hspace*{-1em} \sum_{j=v(A)+1}^{v(B)} \hspace*{-1em} \alpha (v(A)+j)^{\alpha-1}m_B(b_j)\\
~&\geq~ \sum_{j=1}^{v(A)} \alpha (2j)^{\alpha-1}m_A(a_j) + \sum_{j=1}^{v(B)}\alpha  (2j)^{\alpha-1}m_B(b_j)
~=~ 2^{\alpha-1}M(A) + 2^{\alpha-1}M(B)\,. \end{aligned}$$

Let $v_1,\dots,v_{v(G)}$ be a list of $V(G-S)=V(G)$ that attains \vspace*{-1ex}
$$M(G-S) ~=~ \sum_{j=1}^{v(G)} \alpha j^{\alpha-1}m_{G-S}(v_j)\,.$$
We have $d_{G}(v_j)\geq d_{G-S}(v_j)$ and hence $\Delta_j:= m_{G-S}(v_j) -m_{G}(v_j) \in [0,m]$ for every \linebreak $j\in \{1,\dots,v(G)\}$. It is \vspace*{-1ex}
$$\sum_{j=1}^{v(G)} \Delta_j ~\leq~ \sum_{j=1}^{v(G)} d_{G}(v_j)-d_{G-S}(v_j) ~=~ 2|S| ~\leq~ 2k\,.$$
By considering the same list as a list of $V(G)$, we obtain \vspace*{-1ex}
$$M(G)~\geq~ \sum_{j=1}^{v(G)} \alpha j^{\alpha-1}m_{G}(v_j) ~=~ \sum_{j=1}^{v(G)} \alpha j^{\alpha-1}\big(m_{G-S}(v_j)-\Delta_j\big) ~=~ M(G-S)-\sum_{j=1}^{v(G)} \alpha j^{\alpha-1}\Delta_j\,.$$
Since $\sum_{j=1}^{v(G)} \Delta_j \leq 2k$ and the weights $\alpha j^{\alpha-1}$ are decreasing, we obtain an upper bound of the weighted sum over the $\Delta_j \in [0,m]$ by setting $\Delta_j = m$ for $j \leq \big\lceil\frac{2k}{m}\big\rceil$ and $\Delta_j = 0$ otherwise. That is, \vspace*{-1ex}
$$M(G) ~\geq~ M(G-S)-\sum_{j=1}^{\big\lceil\frac{2k}{m}\big\rceil} \alpha j^{\alpha-1}m ~\geq~ 2^{\alpha-1}\big(M(A)+M(B)\big)-m\Big(\frac{4k}{m}\Big)^\alpha\,,$$
using the upper bound on the sum of the first $\big\lceil\frac{2k}{m}\big\rceil \leq \frac{4k}{m}$ weights. \hfill $\Box$\\

\textbf{Theorem \newTh{e Theorem}\\}
\textit{For every $\alpha>0$, there is some $k_0$ such that, for all $k\geq k_0$, every graph with minimum degree at least $k+k^{\frac{1}{2}+\alpha}$ contains a $(k+1)$-edge-connected subgraph.}

\textbf{Proof.~} 
Let $\alpha \in \big(0,\frac{1}{2}\big)$,~ $m = k^{\frac{1}{2}+\alpha} \leq k$~ and~ $\mu = \frac{1}{2} m\sqrt{k}^\alpha$.

Let $G$ be a graph with $v(G)>k$ without a $(k+1)$-edge-connected subgraph. We prove by induction on $v(G)$ that $M(G) \geq \mu$.

A sequence of pairwise distinct vertices $w_1,\dots,w_l \in V(G)$ is called \textit{detachable} from $G$, if for each $j\in\{1,\dots,l\}$, $w_j$ has degree at most $k$ in $G-w_1-\dots-w_{j-1}$. We set
$$z(G) ~:=~ \max\big\{ l \in \mathbb{N} ~:~ w_1,\dots,w_l \in V(G) \text{ is detachable from } G\big\}\,.$$
\textbf{Case I:~} $z(G)\geq \sqrt{k}$.

Thus there is a detachable sequence $w_1,\dots,w_{\big\lceil\sqrt{k}\big\rceil} \in V(G)$. For each $j\in\{1,\dots,\big\lceil\sqrt{k}\big\rceil\}$,\linebreak $d_G(w_j) \leq k+(j-1)$ and hence $m_G(w_j) \geq m-(j-1) \geq m-(\big\lceil\sqrt{k}\big\rceil-1) \geq m-\sqrt{k} \geq \frac{2}{3}m$,\linebreak for $k$ large enough. Using the lower bound on the sum of the first weights, this implies, again for $k$ large enough,
$$M(G) ~\geq~ \sum_{j=1}^{\big\lceil\sqrt{k}\big\rceil} \alpha j^{\alpha-1}m_G(v_j) ~\geq~ \frac{2}{3}m \sum_{j=1}^{\big\lceil\sqrt{k}\big\rceil} \alpha j^{\alpha-1} ~\geq~ \frac{2}{3}m(\sqrt{k}^\alpha-1)~\geq~ \frac{1}{2}m\sqrt{k}^\alpha ~=~ \mu\,.$$ 
\textbf{Case II:~} $z(G)< \sqrt{k}$.

Let $w_1,\dots,w_{z(G)} \in V(G)$ be a longest detachable sequence. Hence, all vertices in the graph $G' = G-w_1-\dots-w_{z(G)}$ have degree greater than $k$. Since $G'$ is a non-empty subgraph of $G$, it is not $(k+1)$-edge-connected. Hence, there is a set $S$ of at most $k$ edges such that $G'-S$ decomposes into two non-empty subgraphs $A$ and $B$. We may assume that all edges of $S$ are between $V(A)$ and $V(B)$.

Since $d_{G'}(a)> k$ for all $a\in V(A)$, we have $v(A) > 1$ and
$$(v(A)-1)k ~=~ v(A)k-k ~< \sum_{a\in V(A)} d_{G'}(a)~-|S| ~= \sum_{a\in V(A)} d_{A}(a) ~\leq~ v(A)(v(A)-1)\,,$$
which yields $v(A)>k$, and similarly $v(B)>k$. As it holds for $G$, also its subgraphs $A$ and $B$ do not contain a $(k+1)$-edge-connected subgraph, so we can apply the induction hypothesis to both $A$ and $B$. 
Together with Lemma \ref{e Theorem}b, we obtain
$$M(G') ~\geq~ 2^{\alpha-1}\big(M(A)+M(B)\big) - m\Big(\frac{4k}{m}\Big)^\alpha ~\geq~ 2^\alpha\mu - m\Big(\frac{4k}{m}\Big)^\alpha\,.$$
Finally we apply Lemma \ref{e Theorem}a to get from $G$ to $G'$ by removing the vertices $w_1,\dots,w_{z(G)}$ one after another
$$M(G) ~\geq~ M(G')-z(G)k^\alpha ~\geq~ 2^\alpha \mu - m\Big(\frac{4k}{m}\Big)^\alpha-k^\frac{1}{2} k^\alpha\,.$$
For $k\to \infty$, both $m\big(\frac{4k}{m}\big)^\alpha = m\, o(k)^\alpha = o(\mu)$ and $k^\frac{1}{2} k^\alpha = m = o(\mu)$. Hence, for $k$ large enough, the sum of both numbers is smaller than $(2^\alpha-1)\mu$, that is, $M(G) \geq \mu$.

This finishes the proof: If $G$ has minimum degree at least $k+m$, then $v(G)>k$ and $M(G) = 0$, a contradiction. \hfill $\Box$\\

The condition $\alpha>0$ is sharp, even for bipartite graphs:

\textbf{Theorem \newTh{e Construction}\\}
\textit{For every constant $C$, there exists a $k\in \mathbb{N}$ and a bipartite graph with minimum degree at least $k+C\sqrt{k}$ without a $(k+1)$-edge-connected subgraph.}

\textbf{Proof.~} We may assume that $C\in \mathbb{N}$ and we choose $k$ as a square such that $\frac{1}{2C}\sqrt{k}\in \mathbb{N}$. Further let $k$ be large enough such that $2^{2C^2} C \leq \sqrt{k}$.

We prove for $j=0,\dots,2C^2$ that there is a bipartite graph $G_j$ with classes $A_j$ and $B_j$ and $X_j\subseteq A_j$ with %\vspace*{-1ex}
\begin{enumerate}
	\item[(i)] $G_j$ does not contain a $(k+1)$-edge-connected subgraph. \vspace*{-1ex}
	\item[(ii)] $|X_j| ~=~ 2^{j+1} C\sqrt{k}$. \vspace*{-1ex}
	\item[(iii)] The vertices in $X_j$ have degree at least $k+\frac{j}{2C}\sqrt{k}$. \vspace*{-1ex}
	\item[(iv)] The vertices in $V(G_j)\setminus X_j$ have degree at least $k+C\sqrt{k}$.
\end{enumerate}
The claim follows with $j=2C^2$.

Define $G_0$ as the graph on the disjoint anticliques $X,B,A,B',X'$ with complete bipartite graphs between each two consecutive sets in the list. Thus $G_0$ is bipartite with classes $A_0 := X\cup A\cup X'$ and $B_0 := B\cup B'$. By choosing $|B|=|A|=|B'|=k$ and $|X|=|X'|=C\sqrt{k}$,\linebreak the $2C\sqrt{k}$ many vertices in $X_0 := X\cup X'$ have degree $k$ and all vertices in $V(G_0)\setminus X_0 = B\cup A\cup B'$ have degree at least $k+C\sqrt{k}$. A $(k+1)$-edge-connected subgraph of $G_0$ cannot contain any vertex of degree at most $k$. Hence, it is completely contained in $G_0-X_0$. Here, the vertices in $B\cup B'$ have only degree $k$, and $A$ itself is too small to contain a $(k+1)$-edge-connected subgraph. Thus, $G_0$ satisfies all claimed properties.

Now let $0\leq j<2C^2$. Suppose we have already constructed $G_j,X_j$ and we want to construct $G_{j+1},X_{j+1}$. Note that $|X_j| ~=~ 2^{j+1} C\sqrt{k} \leq 2^{2C^2}C \sqrt{k} \leq k$.

\begin{center}
\vspace*{-1ex}\includegraphics[width = 0.6\textwidth]{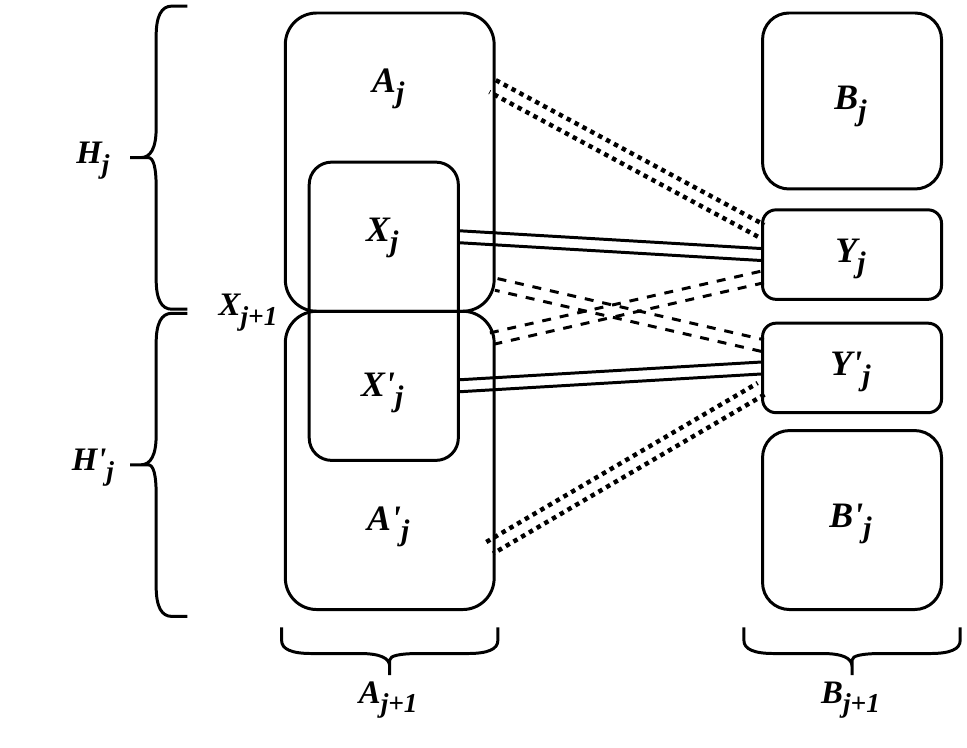}\\
\small Figure 2:~ Construction of $G_{j+1}$ from $G_j$. \vspace*{-1ex}
\end{center}

We first build a graph $H_j$, which is shown in the upper half of Figure 2, from $G_j$ by adding an anticlique $Y_j$ on $\frac{1}{2C}\sqrt{k}\in \mathbb{N}$ vertices. From each vertex in $Y_j$, we add one edge to every vertex in $X_j$ and $k-|X_j|$ more edges to other vertices from $A_j$. In $H_j$, the vertices in $Y_j$ have degree $k$, the vertices in $X_j$ now have degree at least $k+\frac{j+1}{2C}\sqrt{k}$, and all other vertices still have degree at least $k+C\sqrt{k}$. A $(k+1)$-edge-connected subgraph of $H_j$ cannot contain any vertex from $Y_j$, but also $G_j = H_j-Y_j$ does not contain one by the assumption.

We now build $G_{j+1}$ from $H_j$ and one disjoint copy $H'_j$, which is shown in the lower half of Figure 2. Let $A'_j, B'_j, X'_j, Y'_j$ be the copies of $A_j, B_j, X_j, Y_j$ in $H'_j$. From every vertex in $Y_j$, we add $C\sqrt{k}$ edges to vertices in $A'_j$. Similarly, from every vertex in $Y'_j$, we add $C\sqrt{k}$ edges to vertices in $A_j$. Doing so guarantees that $G_{j+1}$ is also bipartite with classes $A_{j+1}=A_j\cup A'_j$ and $B_{j+1}=B_j \cup Y_j \cup B'_j \cup Y'_j$. All together, we have added $2\frac{1}{2C}\sqrt{k}C\sqrt{k} = k$ edges between $H_j$ and $H'_j$. Hence, a $(k+1)$-edge-connected subgraph of $G_{j+1}$ must be completely contained either in $H_j$ or in $H'_j$, which we have precluded.

We set $X_{j+1}:=X_j\cup X'_j$ of size $2|X_j| = 2^{j+1}C\sqrt{k}$. The vertices in $X_{j+1}$ still have degree at least $k+\frac{j+1}{2C}\sqrt{k}$ as in $H_j$. All other vertices, including those in $Y_j\cup Y'_j$, have degree at least $k+C\sqrt{k}$. This finishes the construction. \hfill $\Box$

For the proof, we need $C = \Theta(\sqrt{\log k})$, which might be sharp: For some $D$, a minimum degree of at least $k+D\sqrt{k\log k}$ might be sufficient for a $(k+1)$-edge-connected subgraph. \newpage

\subsection*{Application on graph decomposition}

The theorems from the present paper have some applications to other issues. We want to mention decomposition problems of the following type, which were first investigated by Thomassen \cite{Thomassen}:
What is the lowest number $k=k(s,t)$ such that every $k$-(edge-)connected graph $G$ admits a partition of $V(G)$ into $A$ and $B$ such that $G[A]$ is $s$-(edge-)connected and $G[B]$ is $t$-(edge-)connected? 

\textbf{Corollary.}~~ \textit{Let $k_1,k_2\geq 1$. \vspace*{-1ex}
\begin{enumerate}
\item[(i)] If $G$ is a $(k_1+k_2+1)$-connected graph with minimum degree at least $3(k_1+k_2)-1$, then $V(G)$ has a partition into $C_1$ and $C_2$ such that, for $i\in\{1,2\}$, $|C_i|>2k_i$ and $G[C_i]$ is $(k_i+1)$-connected.
\item[(ii)] If $G$ is a triangle-free $(k_1+k_2+1)$-connected graph with minimum degree at least $2(k_1+k_2)$, then $V(G)$ has a partition into $C_1$ and $C_2$ such that, for $i\in\{1,2\}$, $G[C_i]$ is $(k_i+1)$-connected (and hence $|C_i|\geq 2(k+1)$).
\item[(iii)] If $G$ is a $(k_1+k_2+1)$-edge-connected graph with average degree at least $2(k_1+k_2)+2$, then $V(G)$ has a partition into $C_1$ and $C_2$ such that, for $i\in\{1,2\}$, $|C_i|>2k_i$ and $G[C_i]$ is $(k_i+1)$-edge-connected.
\item[(iv)] There is a function $f(r)=\big(1+o(1)\big)r$ such that the following holds:
If $G$ is a $(k_1+k_2+1)$-edge-connected graph with minimum degree at least $f(k_1+k_2)$, then $V(G)$ has a partition into $C_1$ and $C_2$ such that, for $i\in\{1,2\}$, $G[C_i]$ is $(k_i+1)$-edge-connected.
\end{enumerate} }

\textbf{Proof.~} The approach is due to Thomassen \cite{Thomassen} and Stiebitz \cite{Stiebitz}. \vspace*{-1ex}
\begin{itemize}
\item[(i)] $G$ has minimum degree at least $(3k_1-1)+(3k_2-1)+1$. By a result of Stiebitz \cite{Stiebitz}, $V(G)$ has a partition into $A_1$ and $A_2$ such that, for $i\in\{1,2\}$, $G[A_i]$ has minimum degree at least $3k_i-1$. By Theorem \ref{Theorem}, $G[A_i]$ contains a $(k_i+1)$-connected subgraph on some vertex set $B_i$ of size greater than $2k_i$, for $i\in\{1,2\}$. We choose maximal disjoint $C_1,C_2 \subseteq V(G)$ such that $C_i\supseteq B_i$ and $G[C_i]$ is $(k_i+1)$-connected, for $i\in\{1,2\}$.

Suppose that $X= V(G) \setminus (C_1\cup C_2)\neq \emptyset$. 
Thus $G[C_1\cup X]$ is not $(k_1+1)$-connected, that is, there is some $S_1\subseteq C_1\cup X$ of at most $k_1$ vertices such that $G[(C_1\cup X)\setminus S_1]$ has more than one connected component. One of them has to be disjoint from $C_1$. Let $Y\subseteq X$ be its vertex set.

Now $G[C_2\cup Y]$ is not $(k_2+1)$-connected, that is, there is some $S_2\subseteq C_2\cup Y$ of at most $k_2$ vertices such that $G[(C_2\cup Y)\setminus S_2]$ has more than one connected component. One of them has to be disjoint from $C_2$. Let $Z \subseteq Y$ be its vertex set.
The neighborhood of $Z$ is contained in $S_1\cup S_2$ of size at most $k_1+k_2$, which contradicts the assumption that  $G$ is $(k_1+k_2+1)$-connected.

\item[(ii)] By a result of Kaneko \cite{Kaneko} on triangle-free graphs, $V(G)$ has a partition into $A_1$ and $A_2$ such that, for $i\in\{1,2\}$, $G[A_i]$ has minimum degree at least $2k_i$. Now use Theorem \ref{triangle-free} to obtain disjoint $(k_1+1)$- and $(k_2+1)$-connected subgraphs. Now proceed as before.

\item[(iii)] By a result of Wang and Wu \cite{WW}, $V(G)$ has a partition into $A_1$ and $A_2$ such that, for $i\in\{1,2\}$, $G[A_i]$ has average degree at least $2k_i$. Now use Theorem \ref{e Average} to obtain disjoint $(k_1+1)$- and $(k_2+1)$-connected subgraphs with vertex number greater than $2k_1$ and $2k_2$, respectively. Since Thomassen's argument also works for edge-connectivity, we may proceed as before. 

\item[(iv)] By Theorem \ref{e Theorem} and Theorem \ref{e Average} for small $k$, there exists a function $g(k)=\big(1+o(1)\big)k$,\linebreak $g(k)\leq 2k$ such that for every $k$, every graph with minimum degree at least $g(k)$ contains a $(k+1)$-edge-connected subgraph. We choose 
$$f(r) = \max\big\{g(r-k)+g(k)+1\,\big|\,1\leq k\leq \tfrac{r}{2}\big\}\,.$$
We have that $f(r) = \big(1+o(1)\big) r$ for $r\to \infty$:

If $k\geq \sqrt{r}$, then $g(r-k)+g(k) = \big(1+o(1)\big)(r-k)+\big(1+o(1)\big)k = \big(1+o(1)\big)r$.\linebreak
Otherwise, $g(r-k)+g(k) < \big(1+o(1)\big)(r-k) + 2k \leq \big(1+o(1)\big)r + o(r) = \big(1+o(1)\big)r$.

Now let $G$ be a $(k_1+k_2+1)$-edge-connected graph with minimum degree at least \linebreak $f(k_1+k_2)\geq g(k_1)+g(k_2)+1$. Again by Stiebitz \cite{Stiebitz}, $V(G)$ has a partition into $A_1$ and $A_2$ such that, for $i\in\{1,2\}$, $G[A_i]$ has minimum degree at least $g(k_i)$. We obtain disjoint $(k_1+1)$- and $(k_2+1)$-edge-connected subgraphs. Now proceed as before.\hfill $\Box$ \vspace{1ex}
\end{itemize}

The result of Wang and Wu \cite{WW} may also be applied to obtain average degree versions of (i) and (ii).


\begin{thebibliography}{}


\bibitem{Mader} W. Mader (1972). Existenz $n$-fach zusammenhängender Teilgraphen in Graphen
genügend großer Kantendichte. Abhandlungen aus dem math. Seminar der Universität Hamburg, 37, 86-97. \vspace*{-0.5ex}

\bibitem{Mader e} W. Mader (1971). Minimale $n$-fach kantenzusammenhängende Graphen. Mathematische Annalen, 191, 21–28. \vspace*{-0.5ex}

\bibitem{Matula} D. W. Matula (1978).  Subgraph connectivity numbers of a graph. In: Alavi, Y., Lick, D.R. (eds) Theory and Applications of Graphs: Proceedings (pp. 371-383). Lecture Notes in Mathematics, vol 642. Springer, Berlin, Heidelberg. \vspace*{-0.5ex} 

\bibitem{Mader Conj} W. Mader (1979). Connectivity and edge-connectivity in finite graphs. B. Bollobás (Ed.), Surveys in Combinatorics, Cambridge University Press, 66–95. \vspace*{-0.5ex}

\bibitem{Thomassen} C. Thomassen (1983). Graph decomposition with constraints on the connectivity and minimum degree. Journal of Graph Theory, 7, Issue 2, 165-167.  \vspace*{-0.5ex}

\bibitem{Nguyen} T. H. Nguyen (2024). Highly connected subgraphs with large chromatic number. SIAM Journal on Discrete Mathematics, 38(1), 243-260. \vspace*{-0.5ex}

\bibitem{Bonnet} É. Bonnet, C. Feghali, T. Nguyen, A. Scott, P. Seymour, S. Thomassé, N. Trotignon (2025). Graphs without a 3-Connected Subgraph are 4-Colourable. The Electronic Journal of Combinatorics, P1-26. \vspace*{-0.5ex}

\bibitem{Carmesin} J. Carmesin (2020+). Large highly connected subgraphs in graphs with linear average degree. preprint. \texttt{https://arxiv.org/abs/2003.00942} \vspace*{-0.5ex}

\bibitem{Tutte} W. T. Tutte (1961). On the problem of decomposing a graph into $n$ connected factors. Journal of the London Mathematical Society, 1(1), 221-230. \vspace*{-0.5ex}

\bibitem{N-W} C. St. J. A. Nash-Williams (1961). Edge-disjoint spanning trees of finite graphs. Journal of the London Mathematical Society, 1(1), 445-450. \vspace*{-0.5ex}

\bibitem{Stiebitz} M. Stiebitz (1996). Decomposing graphs under degree constraints. Journal of Graph Theory, 23, Issue 3, 321-324. \vspace*{-0.5ex}

\bibitem{Kaneko} A. Kaneko (1998). On decomposition of triangle‐free graphs under degree constraints. Journal of Graph Theory, 27(1), 7-9.  \vspace*{-0.5ex}

\bibitem{WW} Yan Wang, Hehui Wu (2024). Graph partitions under average degree constraint. Journal of Combinatorial Theory Series B, 165, 197–210. 

\end{thebibliography}
\end{document}